\documentstyle[twoside,12pt]{article}
\newtheorem{prop}{Proposition}[section]
\newtheorem{dfn}[prop]{Definition}   
\newtheorem{theo}[prop]{Theorem}
\newtheorem{conj}[prop]{Conjecture}
\newtheorem{rem}[prop]{Remark} 
\newtheorem{coro}[prop]{Corollary} 
\newtheorem{lem}[prop]{Lemma}   
\newtheorem{exam}[prop]{Example}

\title{Isolated non-normal crossings}
\author{A.Libgober\\
\small Department of Mathematics \\
\small University of Illinois at Chicago\\
\small 851 S.Morgan Str. Chicago, Illinois, 60607 \\
\small e-mail: libgober@math.uic.edu \\}

\begin{document}

\date{}

\maketitle

\begin{abstract}{We describe a multivariable polynomial 
invariant for certain class of 
non isolated hypersurface singularities generalizing 
the characteristic polynomial on monodromy.
The starting point is   
an extension of a theorem due to L\^e and K.Saito on
commutativity of the local fundamental groups of certain 
hypersurfaces. The description of multivariable polynomial
invariants is given in terms of the ideals and polytopes of quasiadjunction 
generalizing corresponding data used in the study of the  
homotopy groups of the complements to projective hypersurfaces
and Alexander invariants of plane reducible curves.} 
\end{abstract}
\section{Introduction}

One of the central results in the study of the topology of algebraic 
hypersurfaces near a singular point is the assertion that the complement 
to the link of a singularity 
is fibered over the circle with the connectivity of the 
fiber depending on the dimension of the singular locus of the 
hypersurface (cf. \cite{milnor}). More precisely, if $X \subset {\bf C}^{n+1}$ 
is a germ of such a hypersurface,
${\rm Sing}X$ is the singular locus of $X$ which we shall assume 
contains the origin $0$, $k={\rm dim}X$ and 
$B_{\epsilon}$ is ball in ${\bf C}^{n+1}$ of a 
sufficiently small radius $\epsilon$ centered at $0$, then 
there is a locally trivial fibration $\partial B_{\epsilon}-
\partial B_{\epsilon} \cap X \rightarrow S^1$. Moreover, the fiber $F_X$
of this fibration is $n-k-1$-connected (i.e $\pi_i(F_X)=0$ for $1 \le  i <n-k$). 
Existence of a locally trivial fibration allows 
to restate this Milnor's result as follows: 
the cyclic cover of $\partial B_{\epsilon}-
\partial B_{\epsilon} \cap X$ corresponding to the kernel of the map 
$\pi_1(\partial B_{\epsilon}-\partial B_{\epsilon} \cap X) 
\rightarrow {\bf Z}$ sending a loop to its linking number with 
$ \partial B_{\epsilon} \cap X$ is $n-k-1$-connected.  In fact this 
cyclic cover is homotopy equivalent to $F_X$. 

In present papers we consider the {\it abelian} rather than cyclic 
covers of the 
complement $\partial B_{\epsilon}-
\partial B_{\epsilon} \cap X$ and study their connectivity.
It is easy to see (cf. lemma \ref{homology}) that $H_1(\partial B_{\epsilon}-
\partial B_{\epsilon} \cap X)={\bf Z}^r$ where $r$ is the number 
of irreducible components of $X$ i.e. a  
hypersurface must be reducible in order that its complement 
will admit abelian but non cyclic covers. 
This forces singularities already in codimension $1$. We shall 
assume that these singularities, inherently present in 
the cases when non cyclic covers do exist, are normal crossings 
except perhaps for the origin. In this case the Milnor fiber has    
the fundamental group isomorphic to ${\bf Z}^{r-1}$ (cf. \cite{LeSaito}
and \ref{comparison}) 
but the universal abelian cover is $(n-1)$-connected i.e. making this 
situation analogous to the case of isolated singularities in cyclic 
case. This result is proven in section \ref{vanishing}.
Standard arguments with hyperplane sections show that if a germ 
of a singularity has 
non-normal crossings only in dimension $k$ then 
$\pi_i(\partial B_{\epsilon}-
\partial B_{\epsilon} \cap X)=0$ for $2 \le i \le n-k-1$. 
The fundamental group $\pi_1(\partial B_{\epsilon}-
\partial B_{\epsilon} \cap X)$ is abelian as was shown by L\^e Dung Trang
and K.Saito (cf. \cite{LeSaito}) and hence is isomorphic to ${\bf Z}^r$.
Therefore section \ref{vanishing} can be viewed as a generalization of 
the result of L\^e-Saito which takes into account dimension of 
the locus of non normal crossings. 

In the case $n=1$ the topology of abelian covers of the complements
to the union of irreducible germs was studied in \cite{alexhodge}.
This case is similar to the one considered in this paper but 
the role of the homotopy group is played by the 
homology of the universal abelian cover. Indeed in the case 
of isolated non normal crossings of a union of germs of $r$ 
hypersurfaces (INNC) the homology and the homotopy of the universal abelian 
covers are isomorphic (cf. section \ref{kgcomplex}). In \cite{alexhodge}
we described an algebro-geometric way to calculate the support of 
the homology of universal abelian cover endowed with the 
structure of ${\bf Z}^r$-module. In this paper we generalize 
this method to the case of isolated non normal crossings. 
More precisely we consider the ${\bf C}[{\bf Z}^r]$-module 
$\pi_n(\partial B_{\epsilon}-
\partial B_{\epsilon} \cap X) \otimes_{\bf Z} {\bf C}$. Its support, 
can be defined as the set of points 
$\chi \in {\rm Spec}{\bf C}[{\bf Z}^r]={{\bf C}^*}^r$
such that: 

$$(\pi_n(\partial B_{\epsilon}-
\partial B_{\epsilon} \cap X) \otimes_{\bf Z} {\bf C}) 
\otimes_{\pi_1(\partial B_{\epsilon}-
\partial B_{\epsilon} \cap X)} {\bf C}_{\chi} \ne 0$$
(we identify ${\rm Spec}[\pi_1(\partial B_{\epsilon}-
\partial B_{\epsilon} \cap X)]$ and 
${\rm Char}(\pi_1(\partial B_{\epsilon}-
\partial B_{\epsilon} \cap X)$
using isomorphism of their rings of regular function 
obtained by viewing elements of 
${\bf Z}^r$ as functions on ${\rm Char}({\bf Z}^r)$;
here ${\bf C}_{\chi}$ is ${\bf C}$ with the action 
of $\gamma \in \pi_1(\partial B_{\epsilon}-
\partial B_{\epsilon} \cap X)$ given by multiplication 
by $\chi(\gamma)$). Notice that the study of the homotopy groups 
of the complements with module structure when the 
fundamental group is ${\bf Z}$ started in 
\cite{annals}, \cite{position}, \cite{cheniot}, \cite{dimca} and 
in the context of non isolated singularities in \cite{tibar}. 

Support loci for homology groups of 
the universal abelian cover of a CW complex may have 
arbitrary codimension in ${\rm Spec}{\bf C}[{\bf Z}^r]$
(cf. e.g. calculations in \cite{charvar}). We show, however, 
that in the case of isolated non normal crossings 
${codim \rm Supp} \pi_n(\partial B_{\epsilon}-
\partial B_{\epsilon} \cap X) \otimes_{\bf Z} {\bf C}$ is 
equal to one.
In fact we calculate several components of the support
of the homotopy group in terms of a resolution of the non normal 
crossing singularity. 
These components correspond to faces of polytopes constructed
here from resolutions and are called {\it the principal
components}. 
This calculation is given in terms of certain 
ideals associated with the singularity generalizing 
the ideals of quasiadjunction studies in 
\cite{alexhodge} in multivariable version and 
in \cite{arcata1}, \cite{loeser} in the case
of one variable polynomial invariants.  

This suggests a polynomial invariant $P(t_1,t_1^{-1},...,t_r,t_r^{-1},X)$
of INNC's defined as 
a generator of the divisorial hull of the first Fitting ideal
of $\pi_n(\partial B_{\epsilon}-
\partial B_{\epsilon} \cap X) \otimes_{\bf Z} {\bf C}$
(it is well defined up to a unit of the ring of Laurent polynomial).
In cyclic case this polynomial is the characteristic 
polynomial of the monodromy.
An interesting problem is to isolate the cases then the 
support of $\pi_n(\partial B_{\epsilon}-
\partial B_{\epsilon} \cap X)$
consist entirely of the principal components described in this paper.

The content of the paper is the following. In the next section we discuss
technical results used in the paper: $(\Gamma,n)$-complexes, 
Leray-Mayer-Vietoris sequence and other preliminary results.
In section \ref{vanishing} we prove mentioned above extension 
of the theorem of L\^e-Saito. Next we show that 
the zero sets of Fitting ideals of homotopy modules, which 
we call the characteristic varieties of the homotopy groups
(cf. \cite{charvar}), can be used to 
calculate the homology of 
abelian covers of $\partial B_{\epsilon}$  
branched over the intersection of 
the latter with a germ of INNC. Characteristic varieties also 
determine the homology of unbranched covers of the complement to 
a germ of INNC in $\partial B_{\epsilon}$ as well as homology 
of the Milnor fiber. 
These results generalize corresponding results in dimension $n=1$ 
which appear in \cite{topandappl} in unbranched and in \cite{sakuma}
in unbranched case. In section \ref{polynomialINNC} we calculate 
what we call ``the principal component''
of the support of the characteristic variety in question in terms
of ideals and polytopes of quasiadjunciton.
In the last section we provide examples and several concluding remarks
relating introduced here ideals and polytopes of quasiadjunction 
to multiplier ideals and log-canonical thresholds. 
Complete calculation of characteristic varieties
remains an interesting problem.

I want to thank A.Dimca for careful reading of 
an earlier version of this paper and his useful comments. I also thank 
the organizers of the VII'th workshop on 
Singularities for their warm hospitality during my visit to San Carlos.

\section{Preliminaries.}

\subsection{Complements to reducible germs}\label{basics}

We continue to use the notions from the Introduction i.e. 
$X=\bigcup_{i=1}^{i=r} D_i \subset {\bf C}^{n+1}$ is a union of germs of
singularities and $B_{\epsilon}$ is a small ball about the origin.
We shall start with the following:

\begin{lem}\label{homology}$H_1(\partial B_{\epsilon}-
\partial B_{\epsilon} \cap X)={\bf Z}^r$
\end{lem}

\noindent This is immediate consequence of the Alexander duality and 
Mayer Vietoris exact sequence since $\partial B_{\epsilon} \cap X$
is a union of $r$ $(2n-1)$-dimensional manifolds intersecting
{\it pairwise} along $(2n-3)$-submanifolds.

Notice that the identification with ${\bf Z}^r=\oplus_{i=1}^r{\bf Z}e_i$ 
depends on the ordering
of irreducible components of $X$. If such selection is made, the isomorphism
in \ref{homology}  sends $e_i$ to a loop having the linking number with 
$i$-th component equal to 1 and zero linking number with the remaining ones. 

Our main result on the homotopy groups $\pi_i (\partial B_{\epsilon}-
\partial B_{\epsilon} \cap \bigcup D_i)$, which we prove in
the next section, is the following: 

\begin{theo}\label{main}
 Let $X=\bigcup_{i=1}^{i=r} D_i \subset {\bf C}^{n+1}$ 
be is union of 
$r$ irreducible germs with normal 
crossings outside of the origin. If $n \ge 2$,  
then $$\pi_1((\partial B_{\epsilon}-
\partial B_{\epsilon} \cap X)={\bf Z}^r \ \
\pi_k((\partial B_{\epsilon}-
\partial B_{\epsilon} \cap X)=0 \ \ \ {\rm for} \ \ 2 \le k <n$$ 
\end{theo}

The term ''normal crossings outside of the origin'' is used in the sense that 
each components is non singular outside of the origin
and intersection of components at each point outside of the 
origin is the normal crossing.
The claim about the fundamental group is immediate consequence 
of the result of L\^e and K.Saito (cf. \cite{LeSaito}) coupled with 
lemma \ref{homology}. 

Here we just point out that theorem \ref{main} yields the following:

\begin{coro}\label{reduction} Let $NN(X)$ be the non-normal locus of 
$X=\bigcup_{i=1}^{i=r} D_i$
i.e. subset of $X$ of points at which $X$ fails to be 
a normal crossing divisor. If codimension of $NN(X)$ in $X$ is greater
than $1$ then:
$$\pi_1((\partial B_{\epsilon}-
\partial B_{\epsilon} \cap X)={\bf Z}^r \ \ {\rm and} \ \ 
\pi_k (\partial B_{\epsilon}-
\partial B_{\epsilon} \cap X)=0 \ \ {\rm for} \ \ 
2 \le k < n-{\rm dim} NN(X)$$ 
\end{coro}

Let us consider generic linear subspace $L$ in ${\bf C}^{n+1}$
having codimension equal to ${\rm dim}NN(X)$ and passing 
through the origin. Then by a theorem of L\^e and Hamm (cf. \cite{LeHamm})
$\pi_k(L \cap \partial B_{\epsilon}-
\partial B_{\epsilon} \cap X)=\pi_k (\partial B_{\epsilon}-
\partial B_{\epsilon} \cap X)$ for $k <{\rm dim} L-1$.
Hence the corollary follows from theorem \ref{main}.

\subsection{Action of $\pi_1(X)$ on higher homotopy groups}
\label{action}
Fundamental group of a CW-complex acts on the homotopy groups
$\pi_i(X)$ via Whitehead product (cf. \cite{White}). 
If $n \ge 2$ and one has
vanishing: $\pi_i(X)=0$ for $2 \le i <n$ then such action 
on $\pi_n(X)$ coincides with the action of the group of 
covering transformation on $H_n(\tilde X)=\pi_n(\tilde X)=\pi_n(X)$
where $\tilde X$ is the universal cover of $X$ (the first equality is
the Hurewicz isomorphism). 

We shall use the functoriality of this action. In particular
if $p: Y \rightarrow X$ is a circle fibration, $\gamma \in 
\pi_1(Y)$ is an element corresponding to the fiber 
and $x \in \pi_j(Y)$ we have: 
$p_*(\gamma \cdot x)=p_*(\gamma) \dot p_*(x)=p_*(x)$.
Since $p_*: \pi_j(Y) \rightarrow \pi_j(X)$ is 
injective for $j \ge 2$ the action of $\gamma$ is trivial.

\subsection{Complexes with vanishing low dimensional homotopy}\label{kgcomplex}

\begin{dfn}(\cite{Dyer}) {A complex $K$ is called a $(\Gamma,k)$
 complex if:
\par 1.$dim K=k$
\par 2. $\pi_1(K)=\Gamma$
\par 3. $\pi_i (K)=0$  for $2 \le i \le k-1$.}
\label{k-complex}
\end{dfn}

For certain $\Gamma$ the homotopy type of a $(\Gamma,k)$ complex is 
determined by the Euler characteristic. Class of such groups includes 
${\bf Z}^r$ (cf. \cite{Dyer}, \cite{Swan}). In particular a $({\bf Z}^r,n)$
complex is homotopy equivalent to a wedge of the $n$-skeleton of 
a $n$-torus $(S^1)^r$ wedged with several copies of $S^n$.

\subsection{Leray-Mayer-Vietoris spectral sequence}
\label{Leray}

Let $X$ be a paracompact topological space 
 and let ${\cal M}=\{ M_i \vert i \in {\cal S}\}$ 
be a locally finite cover of $X$ by closed subsets.
For a sheaf $\cal A$ on $X$ one has the Leray spectral sequence 
(cf. \cite{Godement} Th. 5.2.4) which is  
one of the two spectral sequence corresponding to 
the double complex: $C^{*,*}(X,{\cal M},{\cal A})
=\oplus_{{\rm Card} S=p+1} C^q(\cap_{i \in S} M_i,{\cal A})$
with Cech and usual cochains differentials.
It has as abutment $H^{p+q}(X,{\cal A}) $ and  
$E_2^{p,q}=H^p({\cal M},{\cal H}^q({\cal A}))$.
Here ${\cal H}^q({\cal A})$ is the system of coefficients on 
the nerve of the cover ${\cal M}$ given by 
$S \rightarrow H^q(M_S,{\cal A})$ where $S$ is a subset of the set ${\cal S}$.
Moreover $E_1^{p,q}=C^p({\cal M},{\cal H}^q)=
\oplus_{{\rm Card} S=p+1} H^q(M_S,{\cal A})$.
Applying this to the case of finite union of closed subsets
$X=\bigcup_{i \in {\cal S}} X_i$, a constant sheaf $\cal A$ and denoting 
$X^{[p]}=\coprod X_{i_0} \cap .... \cap X_{i_p}$ we obtain the 
spectral sequence: 

\begin{equation}\label{MV} 
E_1^{p,q}=H^q(X^{[p]}) \Rightarrow H^{p+q}(X) 
\end{equation}

\begin{exam}
In the case ${\rm Card}{\cal S}=2$ this spectral sequence is 
equivalent to the Mayer-Vietoris exact sequence. Indeed, 
only terms $E_1^{0,q}=H^q(X_1) \oplus H^q(X_2)$ and 
$E_1^{1,q}=H^q(X_1 \cap X_2) $ in $E_1$-term of the  
spectral sequence (\ref{MV}) are non trivial.
Moreover $E_2^{0,q}={\rm Ker}H^q(X_1) \oplus H^q(X_2) 
\rightarrow H^q(X_1 \cap X_2)$ and $E_2^{1,q}={\rm Coker}
 H^q(X_1) \oplus H^q(X_2) \rightarrow H^q(X_1 \cap X_2)$.
Since the other differentials must be trivial one has:
$0 \rightarrow E_2^{0,n} \rightarrow H^n(X_1 \cap X_2) 
\rightarrow E_2^{1,n-1} \rightarrow 0$. The exactness of this 
sequence is equivalent to the exactness of the Mayer-Vietories.
\end{exam}

\section{Homotopy vanishing.}\label{vanishing} 

In this section we shall prove the theorem \ref{main}.

\bigskip 

\noindent {\it Step 1. The action of $\pi_1(\partial B_{\epsilon}-\cup_iD_i)$ 
on $\pi_j(\partial B_{\epsilon}-\cup_iD_i)$ is trivial if 
$j \le n-1$.} 
\bigskip 

\noindent Let $T_{\delta,l}$ be a sufficiently small tubular of one of the
hypersurfaces $D_l$ in $B_{\epsilon}$.
Let $L$ be a generic hyperplane passing through the
origin which belongs to $T_{\delta,l}$. 
By a Zariski-Lefschetz type theorem (cf. \cite{LeHamm})
we have surjection for $j \le n-1$ and an isomorphism
for $j < n-1$ (the first and last equalities due to the conical 
structure of singularities and are valid for all $i$): 
\begin{equation}
\pi_j(\partial B_{\epsilon}-\cup_iD_i)\cap L)=
\pi_j((B_{\epsilon}-\cup_iD_i)\cap L) 
\rightarrow \pi_j((B_{\epsilon}-\cup_iD_i))=
\pi_j(\partial B_{\epsilon}-\cup_iD_i) 
\end{equation}
This implies that the second map in the composition:
\begin{equation}
\pi_j(\partial (B_{\epsilon}-\cup_iD_i) \cap L) \rightarrow
\pi_j(T_{\delta,l} \cap \partial B_{\epsilon}-\cup_iD_i) \rightarrow 
    \pi_j(\partial B_{\epsilon}-\cup_iD_i)
\end{equation}

\noindent is surjective for $j \le n-1$

On the other hand we have a locally trivial fibration:
\begin{equation} 
T_{\delta,l} \cap (\partial B_{\epsilon}-\cup_iD_i) 
\rightarrow \partial B_{\epsilon} \cap (D_l -D_l \cap_{i \ne l} D_i)
\end{equation}
which yields that the action on 
$\pi_j(T_{\delta,l} \cap (\partial B_{\epsilon}-\cup_iD_i))$
of the element of 
$\pi_1(T_{\delta,l} \cap (\partial B_{\epsilon}-\cup_iD_i))$
corresponding to the loop which the boundary of the 
transversal to $D_l$ disk in $\partial B_{\epsilon}$
is trivial (cf. sect. \ref{action})
Therefore the action of the element of $\pi_1(\partial B_{\epsilon}-
\cup D_i)$ corresponding to the boundary of a small 2-disk 
normal to $D_l \cap \partial B_{\epsilon}$ in 
$\partial B_{\epsilon}$ on $\pi_i(\partial B_{\epsilon}-\cup D_i)$
is trivial for $j \le n-1$. Since this is the case for all $l, (l=1,...,r)$ 
and 
$\pi_1(\partial B_{\epsilon}-
\cup D_i)$ is generated by these loops the claim follows.

\bigskip 

\noindent {\it Step 2. For $j \le n-1$ one has the isomorphism:}
\noindent \begin{equation}
H_j(\partial B_{\epsilon}-\cup_iD_i)=\Lambda^j({\bf Z}^r)
\end{equation}
and for $j=n$ the surjection:
\begin{equation}
H_n(\partial B_{\epsilon}-\cup_iD_i) \rightarrow \Lambda^n({\bf Z}^r)
\end{equation}

\noindent First notice that exact sequence of 
the pair $(\partial B_{\epsilon}, \partial B_{\epsilon}-\cup_i D_i)$
and the duality yields:
\begin{equation}
 H_j(\partial B_{\epsilon}-\cup_i D_i)=H^{2n-j}(\partial B_{\epsilon} \cap \cup_i D_i)
\end{equation}
The spectral sequence 
\begin{equation}
E_1^{p,q}=H^q(D^{[p]}) \Rightarrow H^{p+q}(\cup_i \partial B_{\epsilon} \cap D_i)
\end{equation}
discussed in \ref{Leray} in the case $X_i=D_i \cap \partial B_{\epsilon}$
has non zeros in the term $E_1$ only for $q=0,n-1-p, n-p,2n-1-2p$.
Moreover $E_1^{p,2n-1-p}=\oplus_{i_0<...<i_p}
H^{2n-1-p}(D_{i_0} \cap ....\cap D_{i_p})=
\Lambda^{p+1}({\bf Z}^r)$ as can be seen by induction from 
the Mayer-Vietoris sequence. 
For $p+q \ge n+1$ the terms in $E_1^{p,q}$ are unaffected 
by subsequent differentials and hence:
\begin{equation}
\Lambda^p({\bf Z}^r)=E_1^{p-1,2n+1-2p}=E_{\infty}^{p-1,2n+1-2p}=
H^{2n-p}(\cup_iD_i)=H_p(\partial
B_{\epsilon}-\cup_i D_i) \ (p \le n-1)
\end{equation}
If $r<n$ then $E_1^{p,q}=0$ for $p \ge r$ and 
for $p=n$ it follows that $E_1^{n-1,1}=\Lambda^n({\bf Z}^r)=0$
which proves the claim if the number of components of the divisor
is less than the dimension of the ambient space.
On the other hand the case of arbitrary $r$ can be reduced to this 
one as follows. Notice that adding sufficiently high powers of 
new variables to equations defining $D_i$'s 
yields the  system of $r$ equations defining  
INNC $\cup \tilde D_i$ in ${\bf C}^N$ such that all sections of 
$\cup \tilde D_i$ by linear subspace through the origin 
are diffeomorphic (since by Tougeron's theorem adding 
terms of sufficiently high degree does not change topological type).
By Lefschetz theorem (cf. \cite{LeHamm})
one has surjection $\pi_n(\partial B_{\epsilon}-\cup_iD_i)
\rightarrow \pi_n(\partial B_{\epsilon}-\cup_i \tilde D_i)=\Lambda^n{\bf Z}^r$.
Taking $N >r$ we obtain the claim.  

\bigskip \noindent

\noindent {\it Step 3. End of the proof.} 

\bigskip 

\noindent We shall finish the proof of the theorem \ref{main} by
 induction over $j$ 
in $\pi_j(\partial B_{\epsilon}-\cup D_i)$. For $j=1$ the claim 
is proven in \cite{LeSaito}. Next assuming vanishing of the 
homotopy up to dimension $j-1$ let us consider the five terms 
exact sequence for a space $\cal X$ on which a group $\pi$ 
acts freely ($(\pi_j)_{\pi}$ is the quotient of covariants):

\begin{equation}\label{groupaction}
H_{j+1} ({\cal X}) \rightarrow  H_{j+1} (\pi) \rightarrow (\pi_j)_{\pi}
 \rightarrow H_j({\cal X}) \rightarrow H_j (\pi) \rightarrow 0   
\end{equation}

\noindent In the case when ${\cal X}$ is the universal abelian cover
of $\partial B_{\epsilon}-\partial B_{\epsilon} \cap \bigcup_i D_i $ 
and $\pi=\pi_1(\partial B_{\epsilon}-\partial B_{\epsilon} \cap 
(\bigcup_{1 \le i \le r} D_i))$ it
follows from the step 2 and the isomorphism  $\pi={\bf Z}^r$ 
that the right homomorphism is isomorphism for $j \le n-1$,
while the left is the isomorphism for $j \le n-2$ and surjective for 
$j=n-1$. On the other hand the step 1 yields that  $(\pi_j(X))_{\pi_1}=
\pi_j(X)$.
Hence $\pi_j(X)=0$ as long as $j \le n-1$ which yields the claim
of theorem \ref{main}.

\section{Homotopy module and the homology of abelian covers.}

We shall describe now the characteristic varieties corresponding 
to INNC, which are the invariants of the first non-vanishing higher 
homotopy group of the complement to the link. 
Then we shall use them to calculate the homology of 
the covers associated with INNC and the homology of Milnor 
fibers of INNC's.

\begin{dfn}\label{homotopy} Homotopy module of an INNC is 
$\pi_n(\partial B_{\epsilon}-\partial B_{\epsilon} \cap 
(\bigcup_{1 \le i \le r} D_i)$
considered as the module over 
$\pi_1(\partial B_{\epsilon}-\partial B_{\epsilon} \cap 
(\bigcup_{1 \le i \le r} D_i)$.
\end{dfn}

\noindent 
Let ${\rm Char}({\pi_1(\partial B_{\epsilon}-\partial B_{\epsilon} \cap 
(\bigcup_{1 \le i \le r} D_i))})=
{\rm Hom}(\pi_1(\partial B_{\epsilon}-\partial B_{\epsilon} \cap 
(\bigcup_{1 \le i \le r} D_i)),{\bf C}^*)$ be the group of 
characters of the fundamental group. This group can be identified with
${\rm Spec}{\bf C}[\pi_1(\partial B_{\epsilon}-\partial B_{\epsilon} \cap 
(\bigcup_{1 \le i \le r} D_i))]$ since any  
 $\gamma \in \pi_1(\partial B_{\epsilon}-\partial B_{\epsilon} \cap 
(\bigcup_{1 \le i \le r} D_i))$ defines the function $\chi 
\rightarrow \chi(\gamma)$ on the group of characters and so does any 
linear combination of the elements of the group algebra.

Let $R={\bf C}[\pi_1(\partial B_{\epsilon}-\partial B_{\epsilon} \cap 
(\bigcup_{1 \le i \le r} D_i))]$. Using the identification of 
$\pi_1(\partial B_{\epsilon}-\partial B_{\epsilon} \cap 
(\bigcup_{1 \le i \le r} D_i))=
H_1(\partial B_{\epsilon}-\partial B_{\epsilon} \cap 
(\bigcup_{1 \le i \le r} D_i))$ with ${\bf Z}^r$ from \ref{basics}
we obtain identification of $R$ with the ring of Laurent 
polynomials ${\bf C}[t_1,t_1^{-1},...,t_r,t_r^{-1}]$.  
The $k$-th Fitting ideal of the homotopy module \ref{homotopy}
is the ideal generated by the minors 
of order $(n-k+1) \times (n-k+1)$ of the map $\Phi : R^m \rightarrow R^n$ 
such that ${\rm Coker} \Phi=
\pi_n(\partial B_{\epsilon}-\partial B_{\epsilon} \cap 
(\bigcup_{1 \le i \le r} D_i) \otimes_{\bf Z} {\bf C}$ (cf. \cite{charvar}).

\begin{dfn} $k$-th characteristic variety $V_k(X)$
of an isolated non-normal crossing
$X=\cup_{i=1}^{i=r} D_i$ is the subset 
of zeros of the $k$-th Fitting ideal in  
${\rm Spec}{\bf C}[\pi_1(\partial B_{\epsilon}-\partial B_{\epsilon} \cap 
(\bigcup_{1 \le i \le r} D_i))]$.
Alternatively (cf. \cite{charvar}) this is the set of $P$
such that 
$${\rm rk}_{\bf C} \pi_n(\partial B_{\epsilon}-
\partial B_{\epsilon} \cap X)\otimes_{\bf Z} {\bf C} 
\otimes_{\pi_1(\partial B_{\epsilon}-
\partial B_{\epsilon} \cap X)} {\bf C}_{P} \ge k$$
with the structure of $\pi_1(\partial B_{\epsilon}-
\partial B_{\epsilon} \cap X)$-module on 
 ${\bf C}_P$ defined as the standard structure on 
a quotient of ${\bf C}[\pi_1(\partial B_{\epsilon}-
\partial B_{\epsilon} \cap X)]$-module (cf. \cite{charvar}).
\end{dfn}

\noindent Notice that if $X_i$ is obtained from $X$ by deleting 
the component $D_i$ we have surjections
$\psi_i: \pi_1(\partial B_{\epsilon}-
\partial B_{\epsilon} \cap X) \rightarrow
\pi_1(\partial B_{\epsilon}-
\partial B_{\epsilon} \cap X_i)$ and hence embedding: 
 $\psi_i^*: {\rm  Spec}[\pi_1(\partial B_{\epsilon}-
\partial B_{\epsilon} \cap X_i)]
\rightarrow {\rm Spec}[\pi_1(\partial B_{\epsilon}-
\partial B_{\epsilon} \cap X)]$.

\begin{dfn} The nonessential part of $k$-th characteristic variety of 
INNC $X$ is the subset of this characteristic variety which 
consist of points which belong to 
${\rm Im}\psi_i^*(V_k(X_i))$.
Complement to the non essential part is 
called the essential part of characteristic variety.
\end{dfn} 

\noindent Now we shall describe the homology of branched and unbranched 
covers in terms of characteristic varieties. For each 
$P \in {\rm Spec}{\bf C}[\pi_1(\partial B_{\epsilon}-
\partial B_{\epsilon} \cap X)]$ let: 

\begin{equation}
   f(P,X)=\{ {\rm max} \  k  \ \vert P \in V_k(X) \}
\end{equation}

\noindent We have the following:

\begin{prop} \label{unbranched}
Let $U_{m_1,...,m_r}$ be unbranched cover 
of $\partial B_{\epsilon}-\partial B_{\epsilon} \cap 
(\bigcup_{1 \le i \le r} D_i)$ corresponding 
to the homomorphism 
$\pi_1(\partial B_{\epsilon}-\partial B_{\epsilon} \cap 
(\bigcup_{1 \le i \le r} D_i))=
{\bf Z}^r \rightarrow \oplus_{1 \le i \le r} {\bf Z}/m_i{\bf Z}$. 
Then $$ {\rm rk}H_p(U_{m_1,...,m_r},{\bf C})=\Lambda^p({\bf Z}^r)
\\\ p \le n-1$$
$${\rm rk}H_n(U_{m_1,...,m_r},{\bf C})
=\sum_{(...,\omega_j,...),\omega_j^{m_j}=1} f((...,\omega_j,...),
\pi_n(\partial B_{\epsilon}-\partial B_{\epsilon} \cap 
(\bigcup_{1 \le i \le r} D_i))$$

\end{prop}

\noindent A proof, similar to the one given in \cite{topandappl}
in the case of links in $S^3$, can be obtained by applying 
the sequence (\ref{groupaction}) to 
the case when $\cal X$ is the universal abelian cover of 
$\partial B_{\epsilon}-\partial B_{\epsilon} \cap 
(\bigcup_{1 \le i \le r} D_i)$ and when $\pi$ is the 
kernel of the homomorphism $\pi_1(\partial B_{\epsilon}-\partial B_{\epsilon} \cap 
(\bigcup_{1 \le i \le r} D_i))=
{\bf Z}^r \rightarrow \oplus_{1 \le i \le r} {\bf Z}/m_i{\bf Z}$.

\bigskip
\noindent In branched case we have:

\begin{prop}\label{branchedcovers} 
Let $V_{m_1,...,m_r}$ be branched cover 
of $\partial B_{\epsilon}$ branched over \newline 
$\partial B_{\epsilon} \cap 
(\bigcup_{1 \le i \le r} D_i)$ with the Galois group 
$G=\oplus_{1 \le i \le r} {\bf Z}/m_i{\bf Z}$. 
For each 
$\chi \in {\rm Char} G$ let $I_{\chi}=\{i \vert 1 \le i \le r, \chi({\bf Z}_{m_i}) \ne 1 \}$ 
where ${\bf Z}_{m_i}$ is the $i$-th summand of $G$. 
Such $\chi$ can also be considered as a character of 
$\pi_1(\partial B_{\epsilon}-\partial B_{\epsilon} \cap 
(\bigcup_{i \in I_{\chi}} D_i))$ 
in which case it will be called reduced and denoted $\chi_{red}$.
Let $V_{\chi}$ be the branched cover of $\partial B_{\epsilon}$
branched over $\partial B_{\epsilon} \cap 
(\bigcup_{i \in I_{\chi}} D_i)$ and having ${\rm Im \chi}=G/{\rm Ker} \chi$
as its Galois group.
Then $$ \pi_p(V_{m_1,...,m_r})=0 \\\ 1 \le p \le n-1$$
$${\rm rk}H_n(V_{m_1,...,m_r},{\bf C})
=\sum_{\chi \in {\rm Char} } f(\chi_{red},
\bigcup_{i \in I_{\chi}} D_i)$$

\end{prop}

\noindent {\it Proof.} The abelian cover $V_{m_1,...,m_r}$ is 
the link of an isolated complete intersection 
singularity:
\begin{equation}\label{CI}
z_1^{m_1}=f_1(x_1,...,x_{n}),...,z_r^{m_r}=
f_r(x_1,...,x_{n})
\end{equation}
where $f_i(x_1,...,x_n)=0, 1 \le i \le r$ is an equation of the divisor 
$D_i$. Hence we obtain the first equality in the proposition
 (cf. \cite{Hamm1}).
\par For $\chi \in {\rm Char} G$ and a $G$-vector space $W$ 
let $W^{\chi}=\{w \in W \vert gw=\chi(g)w \}$. We also will 
denote by $U_{\chi}$ the unbranched cover of 
$\partial B_{\epsilon}-\partial B_{\epsilon} \cap 
(\bigcup_{1 \le i \le r} D_i)$ corresponding to the 
subgroup ${\rm Im} \chi \subset G$. 
We have $U_{\chi}=U_{m_1,...,m_r}/{\rm Ker}\chi$ and 
$V_{\chi}=V_{m_1,...,m_r}/{\rm Ker \chi}$. Hence
\begin{equation}
H_p(V_{\chi},{\bf C})=H_p(V_{m_1,...,m_r},{\bf C})_{{\rm Ker} \chi}
\end{equation}
(on the right is the quotient group of covariants: 
$$H_p(V_{m_1,...,m_r},{\bf C})/(1-g)H_p(V_{m_1,...,m_r},{\bf C})$$
with $g \in {\rm Ker} \chi$). 
The latter implies that 
\begin{equation}
H_p(U_{m_1,...,m_r})^{\chi}=
H_p(U_{\chi})^{\chi} \ \ \ 
H_p(V_{m_1,...,m_r})^{\chi}=
H_p(V_{\chi})^{\chi}
\end{equation} 

If a character $\chi$ is such that $\chi (t_i) \ne 1, q \le i \le r$ 
where $t_i$'s are the standard generators of 
$H_1(\partial B_{\epsilon}-\partial B_{\epsilon} \cap 
(\bigcup_{1 \le i \le r} D_i)$ (cf. section \ref{basics}) then 
\begin{equation}\label{brunbr}
H_n(U_{m_1,...,m_r})^{\chi}=H_n(V_{m_1,...,m_r})^{\chi}
\end{equation}
Indeed in the exact sequence of the pair: 
$$H_{n+1}(V_{m_1,...,m_r},U_{m_1,...,m_r}) 
\rightarrow H_n(U_{m_1,...,m_r}) \rightarrow H_n(V_{m_1,...,m_r}) \rightarrow 
H_n(V_{m_1,...,m_r},U_{m_1,...,m_r}) $$
the terms on the ends are isomorphic to the cohomology 
of preimage of branching locus of the cover. All eigenspaces there
have characters trivial on one of $t_i$ which yields (\ref{brunbr}).
Finally for a character $\chi$ such that $\chi(t_{i_1})=...=
\chi(t_{i_k})=1$ we have:
\begin{equation} 
H_n(V_{m_1,...,m_r})^{\chi}=H_n(V_{\chi})^{\chi}=
H_n(V_{m_1,...,\hat m_i,..m_r})^{\chi}=
H_n(U_{m_1,...\hat m_i,..,m_r})^{\chi}) \\\ (i=i_1,...,i_k)
\end{equation}
(the second equality takes place since for selected $\chi$
the space $V_{\chi}$ constructed for $V_{m_1,...,m_r}$
and $V_{m_1,...,\hat m_i,..m_r}$ are the same.
Now the conclusion follows form the proposition \ref{unbranched}

Similar to the above approach can be used to describe the homology of 
the Milnor fiber of the singularity of $X$ in terms 
of the homotopy module.

\begin{prop}\label{comparison}
The homology of the Milnor fiber $F_X$ of the INNC singularity $X$ 
is given by: 
$$ H_p(F_X,{\bf Z})=\Lambda^p ({\bf Z}^{r-1}) \ \ {\rm for} \  1 \le p < n$$
The action of the monodromy of this homology is trivial.
The multiplicity of $\omega \ne 1$ as a root of the 
characteristic polynomial $\Delta_n(X,t)$ 
 of the Milnor's monodromy on $H_n(F_X,{\bf C})$ is equal to: 
$$m_{\omega}=f((...,\omega,...),X)=
{\rm max} \{ i \vert (\omega, ....,\omega) \in V_i(X)
 \subset {\rm Spec} 
{\bf C}[\pi_1 (\partial B_{\epsilon}-\partial B_{\epsilon} \cap X)] \}$$
\end{prop}

\noindent {\it Proof}. As follows the from Milnor's fibration theorem 
(cf. \cite{milnor}), the Milnor fiber is homotopy 
equivalent to the cyclic cover of  
$\partial B_{\epsilon}-\partial B_{\epsilon} \cap X$ 
corresponding to the subgroup 
${\bf Z}^{r-1}={\rm Ker} \pi \subset \pi_1(\partial B_{\epsilon}-
\partial B_{\epsilon} \cap X)$
where the homomorphism 
$\pi: \pi_1(\partial B_{\epsilon}-
\partial B_{\epsilon} \cap X) \rightarrow {\bf Z}$
is given by $e_i\rightarrow 1$ with generators $e_i$ 
described in \ref{basics}. The fibration of the classifying spaces
corresponding to the homomorphism $\pi$:
${S^1}^r \rightarrow S^1$ is trivial and hence the action of 
the monodromy yielded by $S^1$ on the cohomology of the classifying space
of ${\rm Ker} \pi$ is trivial as well.
We have $\pi_1(F_X)={\rm Ker} \pi={\bf Z}^{r-1}$.
We can calculate the homology of the Milnor fiber using 
the induction starting with $j=2$ and the exact sequence
(\ref{groupaction}) corresponding to the group acting 
freely on the universal abelian cover 
with the quotient $F_X$. In this case it looks as follows:

\begin{equation}
H_{j+1} (F_X) \rightarrow  
H_{j+1} ({\bf Z}^{r-1}) \rightarrow 
(\pi_j(\partial B_{\epsilon}-\partial B_{\epsilon} \cap X)_{{\bf Z}^{r-1}})
 \rightarrow H_j(F_X) \rightarrow H_j ({\bf Z}^{r-1}) \rightarrow 0   
\end{equation}

As long as $j<n$ the middle term is trivial by vanishing theorem 
of section \ref{vanishing}
and we obtains the first claim of the proposition.
For $j=n$ the second from the left term may affect only the root $t=1$
of $\Delta_n(X,t)$. Other roots of $\Delta(X,t)$ are the elements
of characteristic variety of ${\bf C}[t,t^{-1}]$-module $H_n(F_X,{\bf C})$ 
and the second claim follows since the map
${\rm Spec}{\bf C}[t,t^{-1}] \rightarrow 
  {\rm Spec} 
{\bf C}[\pi_1 (\partial B_{\epsilon}-\partial B_{\epsilon} \cap X)]$
corresponding to homomorphism $\pi$ is given by $\omega \rightarrow
(\omega, ..., \omega)$ (cf. \cite{charvar}).

\begin{prop}\label{essential} 
Let $X$ be a germ of INNC having $r-1$ irreducible
components and let $D$ be a hypersurface 
with isolated singularity at the origin. Let $i_D$ be the 
surjection: $\pi_1 (\partial B_{\epsilon}-\partial B_{\epsilon} \cap 
(X \cup D)) \rightarrow 
\pi_1 (\partial B_{\epsilon}-\partial B_{\epsilon} \cap 
X)$ inducing the map of the groups of characters: 
 $$i_D^*: 
{\rm Char} \pi_1 (\partial B_{\epsilon}-\partial B_{\epsilon} \cap 
X) \rightarrow {\rm Char} \pi_1 (\partial B_{\epsilon}-
\partial B_{\epsilon} \cap (X \cup D))$$
Then $$i_D^*(V_i(X) \cap {\rm Tor \ Char} 
\pi_1 (\partial B_{\epsilon}-\partial B_{\epsilon} \cap 
X) \subset V_i(X \cup D)$$
\end{prop}

\noindent {\it Proof.} First notice that if $Y$ is a link of INCC 
and $h$ is a diffeomorphism of $Y$ 
leaving irreducible components and their orientations invariant then
\begin{equation}\label{connectivity}
           h_* \vert_{H_i(Y,{\bf Z})}=id \ \ {\rm for} 
\  [{{{\rm dim} Y} \over 2}]+1 < i \le {\rm dim} Y  
\end{equation}
(here [...] denotes the integer part). This follows from the  
Mayer-Vietoris spectral sequence considered in section \ref{vanishing}.
Indeed, the action of $h_*$ is trivial on the terms $E_1^{p,q}$ with 
$p+q> [{{{\rm dim} Y} \over 2}]+1$ since each such non zero $E_1^{p,q}$ 
is generated by the fundamental class of a link of isolated 
hypersurface singularity.

The map $i_D^*$, in coordinates dual to coordinates from 
the identification described in \ref{basics}:  
$H_1(\partial B_{\epsilon}-\partial B_{\epsilon} \cap X,{\bf Z})=
{\bf Z}^{r-1}$, takes 
$(\omega_1,...,\omega_{r-1}) \rightarrow (\omega_1,...,\omega_{r-1},1)$
(where $\omega_i$ is a root unity of degree $m_i$).
 Let $\pi_{m_1,...,m_{r-1}}: U_{m_1,...,m_{r-1}} \rightarrow 
\partial B_{\epsilon}-\partial B_{\epsilon} \cap X$ 
be the covering map. 
We have the exact sequence:
$$H_n(U_{m_1,...,m_{r-1}}-\pi^{-1}_{m_1,...,m_{r-1}} (D \cap
(\partial B_{\epsilon}-\partial B_{\epsilon} \cap X))) \rightarrow   
 H_n(U_{m_1,...,m_{r-1}}) \rightarrow$$
\begin{equation}\label{extracomp}
\rightarrow 
 H_n(U_{m_1,...,m_{r-1}}-\pi^{-1}_{m_1,...,m_{r-1}} (D \cap
(\partial B_{\epsilon}-\partial B_{\epsilon} \cap X)),   
 U_{m_1,...,m_{r-1}}) \rightarrow 
\end{equation}

Using Lefschetz duality one obtains that the above relative 
homology group is isomorphic to  
\begin{equation}
H^{n+1}(\bar \pi^{-1}_{m_1,...,m_{r-1}}(D \cap B_{\epsilon}),
\bar \pi^{-1}_{m_1,...,m_{r-1}}(D \cap X \cap B_{\epsilon}))
\end{equation}
$\pi^{-1}_{m_1,...,m_{r-1}}(D \cap \partial B_{\epsilon})$, which is 
the abelian cover of $D \cap \partial B_{\epsilon}$ branched 
over $D \cap X \cap \partial B_{\epsilon}$
can be identified with 
the link of singularity 
$$u_1^{m_1}=g_1(x_1,...,x_{n+1}),...,u_r^{m_r}=g_r(x_1,...,x_{n+1}), 
f(x_1,...,x_{n+1})=0$$
where $g_i$ are the equations of the components of $X$ and $f$ is the 
equation of $D$. Similarly the preimage of the ramification locus 
$\bar \pi_{m_1,..,m_{r-1}}$ can be identified with the link of 
INNC.  Hence the exact sequence of pair and (\ref{connectivity}) 
yield that the action of the covering 
transformation on the last term in  
(\ref{extracomp}) is trivial. Therefore the multiplicity of the 
character corresponding to $(\omega_1,...,\omega_{r-1})$ 
in $H_n(U_{m_1,..,m_{r-1}})$  
does not exceed the multiplicity of $(\omega_1,..,\omega_{r-1},1)$
in $H_n(U_{m_1,..,m_{r-1}}-\pi^{-1}_{m_1,...,m_{r-1}}(D))$. 
Now the result follows from Prop. \ref{unbranched}.

\begin{rem} The proposition \ref{essential} shows that non empty components 
of an INNC $X$ containing the torsion points 
(in particular principal components
described in the next section) always define components of 
characteristic varieties of homotopy modules of $X'$ obtained from 
$X$ by adding extra-components.
\end{rem}  

\section{Polynomial invariants of INNC.}\label{polynomialINNC}

In this section we show how one can calculate the components of characteristic
varieties using resolution of singularities.

Recall first the results on the mixed Hodge structure on 
the cohomology of links of isolated singularities (\cite{Arcata}).
Let $V$ be a germ of a complex $m$-dimensional space having 
isolated singularity at $p \in V$ and $L$ be the link 
of this singularity. Let $\pi: \tilde V \rightarrow V$ be a 
resolution such that the exceptional locus $E$ 
 is a divisor with normal crossings. 
The cohomology groups $H^*(L)$, which can also be viewed as 
the cohomology with support $H^{*+1}_{\{p\}}(V)$,  
carry the canonical mixed Hodge structure (cf. \cite{Arcata} (1.7)). 
Its Hodge numbers we denote 
 $$h^{ipq}(L)={\rm dim}Gr_F^qGr_{\bar F}^pGr_W^{p+q} H^{i}(L)
   = {\rm dim}Gr_F^qGr_{\bar F}^pGr_W^{p+q} H^{i+1}_{\{p\}}(V)$$
We will need the expression of $h^{ipq}$ in terms 
of differentials on the smooth 
part of $V$. We have: 
$$ h^{(m-1)pq}(L)=h^{(m-1)pq}(E) \\\ for \ \  p+q<m-1$$
\begin{equation}
h^{(m-1)p(m-1-p)}(L)=h^{(m-1)p(m-1-p)}(E)-h^{(m+1)(m-p)(p+1)}(E)
\end{equation}
(cf. \cite{Arcata} Cor.(1.12), the first exact sequence and 
$Gr_{p+q}^WH^i_E(\tilde V)=0$ for $p+q \ne i, i \le m$).
Notice that Mayer Vietors exact sequence and induction over the 
number components of $E$ yield that $h^{m+1,i,j}(E)=0$ if 
$i>m-1$ or $j>m-1$ since each irreducible 
component of $E$ has dimension equal to $m-1$. 
Therefore:
\begin{equation}
h^{(m-1)k0}(L)=h^{(m-1)k0}(E) \ \ \ \ \ 0 \le k \le m-1 
\end{equation}
Moreover, the Fujiki duality ${\rm Hom}(H^i(E),{\bf Q}(-m))=
H^{2m-i}_E(\widetilde V)$ (cf. \cite{Arcata} (1.6))
allows to interpret the cohomology of $E$ in terms of 
cohomology of with support on $E$:
 $$h^{(m-1)pq}(E)=h^{(m+1)(m-p)(m-q)}_E(\tilde V)$$ where
$$h^{(m+1)ab}_E={\rm dim} Gr_F^aGr_{\bar F}^bGr_W^{a+b} H^{m+1}_E(\tilde V)$$
Filtrations on $H^{m+1}_E$ have the form:
 $$H^{m+1}_E=F^0 \supseteq  F^{1} \supseteq... \supseteq F^{m+1} \supseteq 0$$
 
 $$ 0=W_{m} \subseteq W_{m+1} ... \subseteq W_{2m}=H^{m+1}_E $$

\noindent Therefore, since $\sum_{p \ge m+2, a+b=w} h^{(m+1)ab}=
{\rm rk} F^{m+2}(W_w(H^{m+1}_E)
/W_{w-1}(H^{m+1}_E))=0$, we have: 
\begin{equation}
F^{m+1}(W_w(H^{m+1}_E)/W_{w-1}(H^{m+1}_E))=
h^{(m+1)(m+1)(w-m-1)}_E   \ \ \ (m+1 \le w \le 2m)
\end{equation}
\noindent Hence:
$$h^{(m+1)(m+1)0}_E+...+h^{(m+1)(m+1)(l-m-1)}_E=
{\rm dim}W_{m+l}F^{m+1}H^{m+1}_E={\rm dim}H^0(W_l\Omega^m_{\widetilde V}({\rm log}E)/\Omega^m_{\widetilde V}) $$

The last equality is a consequence of the bifiltered isomorphism 
$H^*_E(\tilde V)=$ \newline  
${\bf H}^*(\Omega^{*-1}_{\widetilde V}({\rm log} E)/
\Omega^{*-1}_{\widetilde V})$ (i.e. respecting weight and Hodge filtrations)
with the Hodge filtration on the latter
obtained from the spectral sequence  starting with \newline 
$H^q(W_l\Omega^p_{\widetilde V}({\rm log}E)/\Omega^m_{\widetilde V})$
and abuting to ${\bf H}^{*}(W_l\Omega^{*-1}_{\widetilde V}({\rm log} E)/
\Omega^{*-1}_{\widetilde V})$ 
and degeneration of this spectral sequence (cf. \cite{Arcata},(1.6)).
We have $H^1(\Omega^m_{\tilde V})=0$ by a theorem of 
Grauert and Riemenschneider (\cite{GR}).
Hence we obtain:
$$
h^{(m+1)(m+1)l}_E={\rm dim Coker}(H^0(\Omega^m_{\widetilde V}) 
\rightarrow H^0(W_l\Omega^m({\rm log}E)))
-{\rm dim Coker}(H^0(\Omega^m_{\widetilde V}) \rightarrow$$
\begin{equation} \rightarrow  
H^0(W_{l-1}\Omega^m({\rm log}E)))
\end{equation}

We shall apply this identity to calculation of the Hodge numbers 
of branched covers $V_{m_1,...,m_r}$ of $\partial B_{\epsilon}$
(cf. Prop.  \ref{branchedcovers}).

\begin{dfn} Let $f_i=0$ be the equation of divisor $D_i$
and let $\pi: \tilde {\bf C}^{n+1} \rightarrow {\bf C}^{n+1}$ 
be a resolution of the singularities of $\bigcup D_i$ 
(i.e. the proper preimage of the latter is a normal crossings divisor).
Let ${\cal V}_{m_1,...,m_r}$ be the singularity (\ref{CI}) 
having $V_{m_1,...,m_r}$ as its link.  
Let $\tilde V$ be a normalization of 
${\tilde {\bf C}^{n+1}} \times_{{\bf C}^{n+1}} {\cal V}_{m_1,..,m_r}$ 
(cf. (\ref{CI}))
The ideal of quasiadjunction of type $(j_1,...,j_r \vert m_1,...,m_r)$
is the ideal ${\cal A}(j_1,...,j_r \vert m_1,....,m_r)$
of germs $\phi \in {\cal O}_{0,{\bf C}^{n+1}}$ such that $(n+1)$-form:
\begin{equation}\label{diffform}
\omega_{\phi}={{\phi z_1^{j_1} \cdot ...\cdot z_r^{j_r} dx_1 \wedge ...
\wedge dx_{n+1}} \over {z_1^{m_1-1} \cdot ... \cdot z_r^{m_r-1}}}
\end{equation}
on the non singular locus of $V_{m_1,...,m_r}$  after the pull back on 
$\tilde V$ extends over the exceptional set.

The $l$-th ideal of log-quasiadjunction  
${\cal A}_l({\rm log} E)(j_1,...,j_r \vert m_1,....,m_r)$
is the ideal of $\phi \in {\cal O}_{0,{\bf C}^{n+1}}$  such that 
the the pull back of the corresponding form $\omega_{\phi}$ 
on $\tilde V$ is ${\rm log}$-form on $(\tilde V, E)$ 
having weight at most $l$.
\end{dfn}

Collection of ideals of log-quasiadjunction (which is finite, 
cf. \cite{charvar}) allows to define the following 
collection of subsets in the unit cube.

\begin{prop}\label{quasiadjunction}
There exist a collection of subsets ${\cal P}_{\kappa}, (\kappa \in {\cal K})$
in the unit cube $${\cal U}=\{(x_1,....,x_r) \vert 0 \le x_i \le 1 \}$$
in ${\bf R}^r$ and a collection of affine hyperplanes 
${l}_i(x_1,...,x_r)=\alpha_i$ such that 
each ${\cal P}_{\kappa}$ 
is the boundary of the polytope consisting of solutions to
 the system of inequalities:
$${ l}_i(x_1,...,x_r) \ge \alpha_i$$ 
and such that 
\begin{equation}\label{vector}
({{j_1+1} \over {m_1}}, ... ,{{j_r+1} \over {m_r}}) \in {\cal U}
\end{equation} 
belongs to ${\cal P}_{\kappa}$ if and only if
 
\begin{equation}\label{loginequality} 
 {\rm dim} {\cal A}({\rm log} E)(j_1,...,j_r \vert m_1,....,m_r)/
{\cal A}(j_1,...,j_r \vert m_1,....,m_r) \ge 1 
\end{equation} 

Moreover 
\begin{equation}\label{kloginequality} 
 {\rm dim} {\cal A}_{l}({\rm log} E)(j_1,...,j_r \vert m_1,....,m_r)/
{\cal A}_{l-1}(j_1,...,j_r \vert m_1,....,m_r) \ge k 
\end{equation}
if only if (\ref{vector}) belongs to a collection of certain faces 
${\cal P}_{\kappa,\iota}^{k,l} (\iota \in {\cal I}^{k,l})$
of polytopes ${\cal P}_{\kappa}$.

\end{prop}

The proof is completely analogous to the proof of the proposition 2.2
in \cite{alexhodge} or proposition 2.3.1 in \cite{charvar} 
and details 
will be omitted. We shall record, however, the equations for the affine 
hyperplanes mentioned in proposition \ref{quasiadjunction}. 
Let $\pi: \tilde {\bf C}^{n+1} \rightarrow {\bf C}^{n+1}$
be a resolution of singularities of $\bigcup D_i$ i.e. the irreducible 
components $E_j$ of the exceptional set of $\pi$ and  
the proper $\pi$-transforms of $D_i$'s form a normal crossing divisor.
Let, as in (\ref{CI}), $f_i$ be the defining equations of $D_i$.
Let ${\rm ord}_{E_j} \pi^*(f_i)=a_{i,j}$ and 
${\rm ord}_{E_j} \pi^*(dx_1 \wedge ... \wedge dx_{n+1})=c_j$.
Finally for $$\phi \in {\cal A}({\rm log} E)(j_1,...,j_r \vert m_1,....,m_r)/
{\cal A}(j_1,...,j_r \vert m_1,....,m_r)$$ for some 
array $(j_1,...,j_r \vert m_1,...,m_r)$ let $e_j(\phi)={\rm ord}_{E_j}
\pi^*(\phi)$.
Then (cf. \cite{alexhodge}):
\begin{equation}\label{faceequation}
     a_{1,j}(1-x_1)+...a_{r,j}(1-x_r)=e_j(\phi)+c_j+1 
\end{equation} 
is equation of the hyperplane containing a face of a polytope of quasiadjunction.
Moreover, $\phi \in {\cal A}({\rm log} E)(j_1,...,j_r \vert m_1,....,m_r)$
if and only if:
\begin{equation}\label{faceinequality}
a_{1,j}(1-{{j_1} \over {m_1}})+...+
a_{r,j}(1-{{j_r} \over {m_r}}) \le e_j(\phi)+c_j+1
\end{equation}
In the case of ideal of quasiadjunction we have:
$\phi \in {\cal A}(j_1,...,j_r \vert m_1,....,m_r)$
if and only if: 
\begin{equation} 
a_{1,j}(1-{{j_1} \over {m_1}})+...+
a_{r,j}(1-{{j_r} \over {m_r}}) < e_j(\phi)+c_j+1
\end{equation}

\noindent Hence the existence of $\phi$ satisfying (\ref{loginequality}) 
is equivalent to  $({{j_1+1} \over {m_1}}, ... ,{{j_r+1} \over {m_r}})$ 
being solution to all inequalities (\ref{faceinequality}) with 
at least one of them 
being solution to equality (\ref{faceequation}). 
Moreover, for $\phi \in {\cal A}({\rm log} E)(j_1,...,j_r \vert m_1,....,m_r)$
the corresponding form $\omega_{\phi} \in \Omega_{\widetilde V}(log E)$ 
has weight at most $l$ 
if and only if $({{j_1+1} \over {m_1}}, ... ,{{j_r+1} \over {m_r}})$
satisfies at most $l$ equalities corresponding to components of exceptional 
divisor having $l$-fold intersections. This yields the second claim 
of the propositon.

\begin{dfn} The polytopes ${\cal P}_{\kappa}$
 existence of which is asserted in the 
Proposition \ref{quasiadjunction} are called the 
polytopes of quasiadjunction.
\end{dfn}

\noindent One has the following description of the smallest element 
in the Hodge filtration of the cohomology of a link of singularity  
(\ref{CI}) in terms of polytopes of quasiadjunciton and their faces.

\begin{theo} A character of 
$\pi_1(\partial B_{\epsilon}-\partial B_{\epsilon} \cap 
(\bigcup_{1 \le i \le r} D_i)$ acting on  
$W_l(F^{n}H^n(V_{m_1,....,m_r}))$ via the action of the 
Galois group has the eigenspace of dimension at least $k$
if and only if it has the form: 
$$(exp 2 \pi \sqrt {-1} a_1,...,exp 2 \pi \sqrt {-1} a_r)$$
where $(a_1,....,a_r)$ belongs to one of the 
faces ${\cal P}^{k,l}_{\kappa,\iota}$ of a polytope ${\cal P}_{\kappa}$ 
of quasiadjunction of $\bigcup D_i$.
\end{theo}

The proof is similar to the proof in the case of 
reducible curves given in \cite{alexhodge} and will be omitted.

\begin{dfn} Zariski closure of the image of a face of quasiadjunction 
${\cal P}_{\kappa,\iota}^{k,l}$ 
will be called a principal component of the characteristic variety $V_k$.
Polynomial invariant of $X$ is a generator of the divisorial 
hull of the first Fitting ideal of 
$\pi_n(\partial B_{\epsilon}-\partial B_{\epsilon} \cap X) \otimes {\bf C}$.
\end{dfn}

For $k>1$ the varieties $V_k$ may have codimension greater than one.
We don't know if the union of principal components is equal to 
the characteristic variety $V_k$ of $\pi_n(\partial B_{\epsilon}-\partial B_{\epsilon} \cap X) \otimes {\bf C}$. However we shall conjecture the following
(i.e. that the situation is similar to the one in the case of curves
(cf. \cite{alexhodge})):

\begin{conj} Characteristic variety is a union of translated subtori 
of \newline
${\rm Spec}{\bf C}[\pi_1(\partial B_{\epsilon}-\partial B_{\epsilon} \cap X)]$ 
with each translations given by a point of finite order. 
\end{conj}

\noindent This is the true in irreducible case since the roots of the characteristic 
polynomial of the monodromy are roots of unity. Conjecture is also true in the case of curves
essentially since algebraic links are iterated torus links (cf. \cite{alexhodge}).

\section{Examples and Remarks}

\begin{exam} Let $L_i(x_0,...,x_n)=0 \ \ i=1,...r$ be generic linear forms.
Then the subset $\bigcup D_i$ in ${\bf C}^{n+1}$ given by:
\begin{equation}\label{generic}
L_1 \cdot ...\cdot L_r=0
\end{equation}
is the divisor with normal crossings except for 
the origin. The complement to divisor (\ref{generic}) is 
${\bf C}^*$-fibration over the complement in ${\bf P}^n$ 
to the divisor $\bigcup P(D_i)$
given by the same equation (\ref{generic}) in homogeneous coordinates of 
${\bf P}^n$. Hence 
\begin{equation}
\pi_n({\bf C}^{n+1}-\bigcup D_i)=
\pi_n({\bf P}^n-\bigcup PD_i)
\end{equation}

Here the left hand side is considered as the module over
${\bf Z}[t_1,t_1^{-1},...,
t_r,t_r^{-1}]$ while the right hand side is
the   ${\bf Z}[t_1,t_1^{-1},...,
t_r,t_r^{-1}]/(t_1...t_r-1)$-module.

The complement ${\bf P}^n-\bigcup PD_i$
latter has the homotopy type of the $n$-skeleton of the torus
$(S^1)^{r-1}$ (cf. \cite{Hattori}) and hence is a $({\bf Z}^r,n)$
complex (cf. section \ref{kgcomplex}). 
The universal abelian cover has cell structure obtained
by removing from ${\bf R}^r$ the ${\bf Z}^r$-orbits of interiors 
of faces of dimension greater than $n$ in the unit cube. This 
shows that the chain complex of the universal abelian 
cover is the truncated Koszul complex 
corresponding to the system of parameters
$(t_-1,....,t_r-1)$ over the ring 
of Laurent polynomials $R={\bf Z}[t_1,t_1^{-1},...,
t_r,t_r^{-1}]/(t_1...t_r-1)$:
\begin{equation}
0 \rightarrow {\Lambda}^n(R^r) \rightarrow {\Lambda}^{n-1}(R^r) 
\rightarrow .... {\Lambda}^1(R) \rightarrow (t_1-1,...,t_r-1)R \rightarrow 0
\end{equation}
This yields that 
\begin{equation} \pi_n(S^{2n+1}-\bigcup D_i)=H_n(S^{2n+1}-\bigcup D_i)=
{\rm Ker} {\Lambda}^n(R^r) \rightarrow {\Lambda}^{n-1}(R^r)
\end{equation}
Hence we 
obtain the presentation for the homotopy group:
\begin{equation}
\Lambda^{n+1}({\bf Z}[t_1,t_1^{-1},...,t_r,t_r^{-1}]/(t_1...,t_r-1)^r)
\rightarrow \Lambda^n({\bf Z}[t_1,t_1^{-1},...,t_r,t_r^{-1}]/(t_1...,t_r-1)^r)
\rightarrow 
\end{equation}
$$\pi_n({\bf C}^{n+1}-\bigcup D_i) \rightarrow 0 $$

In particular the support of $\pi_n ({\bf C}^{n+1}-\bigcup D_i)$ 
is given by $t_1....t_r=1$

\end{exam}

\begin{exam} More generally, consider the cone over union of non 
singular hypersurfaces $D_{i_0},...,D_{i_r}$
in ${\bf P}^n$ having the degrees ${\rm {deg}} D_{i_k}=d_k$
and forming there a divisor with normal crossings. Indeed, ${\bf P}^n$
is $({\bf Z}^r/(d_1,...,d_r),n)$ complex such that 
passing to finite cover yields the $K$-complex over a finite abelian 
cover which by \cite{Swan} is homotopy equivalent to the 
wedge of spheres $S^{n}$. Hence the support of the homotopy 
group is given by 
\begin{equation}\label{charvar1}
t_1^{d_1} \cdot ... \cdot t_r^{d_r}-1=0
\end{equation}

\par This result can be obtained using the method of previous section.
Indeed, the cone singularity can be resolved by single blow up 
which yields the polytope of quasiadjunction with the faces: 
\begin{equation}\label{conecase}
 d_1x_1+...d_rx_r \ge l
\end{equation}
Indeed if $f_i$ are the defining equations of $D_i$ and 
$\phi \in {\cal M}^{r-l-1}$ we obtain
${\rm mult}_{E}\pi^*(f_i)=d_i$, ${\rm mult}_{E}{\phi}=r-l-1$,
${\rm mult}\pi^*{dx_1 \wedge ... \wedge dx_{n+1}}=n$.
It follows that the ideals of quasiadjunction are the powers
of the maximal ideal and hence the claim (\ref{conecase}). 
The exponential map takes the union of the faces of quasiadjunction 
(\ref{conecase})
into the set having as its Zariski closure (\ref{charvar1}).

\begin{rem} {\bf  Multiplier ideals and log-canonical thresholds
in terms of ideals of quasiadjunction.} 
Recall that for a divisor $D$ in $X$
one defines the multiplier ideal as 
${\cal J}(D)=f_*({\cal O}(K_Y-f^*(K_X)-\lfloor E \rfloor))$ where 
$f: Y \rightarrow X$ is an embedded resolution of pair $(X,D)$.
Consider divisor  
$D_{\gamma_1,....,\gamma_r}=f_1^{\gamma_1} \cdot ... \cdot f_r^{\gamma_r}$
on ${\bf C}^n$. 
We have ${\cal J}(D_{\gamma_1,...,\gamma_r})=
{\cal A}(j_1,...,j_r \vert m_1,....,m_r)$ where $\gamma_i=1-{{j_i+1}
\over {m_i}}$

The face of quasiadjunction  corresponding to $\phi=1$ can be described 
as the closure of the set of points $(\gamma_1,...,\gamma_r)$
in the unit cube  such that the divisor $\sum_i (1-\gamma_i)D_i$ is 
log canonical i.e. the face is formed by log-canonical thresholds. 

\end{rem}

\end{exam}

\end{document}